\documentclass[12pt]{article}

\usepackage{amsthm}
\usepackage{amssymb}
\usepackage{amsfonts}

\newtheorem{thm}{Theorem}[section]

\newtheorem{Lemm}[thm]{Lemma}
\newtheorem{defn}[thm]{Definition}

\newtheorem{Exam}[thm]{Example}

\newcommand{\mult}{\mathrm{mult}}
\newcommand{\emb}{\mathrm{emb}}
\newcommand{\ind}{\mathrm{index}}

\begin{document}
\setcounter{section}{-1}
\begin{center}
{\Large On the multiplicity of terminal singularities  
on threefolds}\\
\end{center}
\begin{center}
 Nobuyuki Kakimi
\end{center}

\begin{flushleft}
 Department of Mathematical Sciences, University of Tokyo, Komaba, Meguro, 
 Tokyo 153, Japan ( e-mail:kakimi@318uo.ms.u-tokyo.ac.jp )  
\end{flushleft}

\begin{flushleft}
Abstract. We give the multiplicity of terminal singularities on
threefolds by simple calculation. 
Then we obtain the best inequalities for
the multiplicity and the index. 
By using this, we can improve 
the boundedness number of terminal weak $\mathbb{Q}$-Fano $3$-folds
in [KMMT, Theorem 1.2].
Furthermore, we can extend [K, Theorem 3.6] 
for Fujita freeness conditions to nonhypersurface terminal singularities.  
\end{flushleft}

\section{Introduction}
Our results are the multiplicity of terminal singularities
and the best inequalities 
for the multiplicity and the index of terminal singularities 
on threefolds. 
Our results are partially generalizations of Artin [A]'s result, that, 
for a normal surface $S$, a rational singular point $p$ of $S$,
$\emb {\dim}_p S = {\mult}_{p} S + 1 $.

We shall prove the following results in this paper:
(Theorem 2.1) 
Let $(X,p)$ be a $3$-fold terminal singular point 
over $\mathbb{C}$. Then,
for all integers $k$
${\dim} m_{X_p}^k/m_{X_p}^{k+1} = {\mult}_{p}X \cdot k(k+1)/2 + k+1$
and  
${\emb}{\dim}_{p} X = {\mult}_{p} X+2$
and 
${\mult}_{p} X \leq {\ind}_{p}X  +2$ 
$(\textrm{if }  {\ind}_{p}X = 1, \textrm{then } {\mult}_{p} X  \leq 2)$.

We can improve [KMMT Theorem 1.2 (2)] 
by (Theorem 2.1) to the following:
(Theorem 3.4)
Let $X$ be a terminal weak $\mathbb{Q}$-Fano $3$-fold.
Then the following hold.
$(1)$ $ -K_x \cdot c_2(X) \geq 0$, and hence $I(X)|24!$. 
$(2)$ Assume further that the anti-canonical morphism 
$g:X \rightarrow  \bar{X}$ does not contract any divisors. 
Then $(-K_X)^3 \leq 6^3 \cdot (2 + 24!) $.
$(3)$The terminal $\mathbb{Q}$-Fano $3$-folds are bounded.

We also can extend [K 3.6] by (Theorem 2.1) 
to the following:
(Theorem 4.1)
Let $X$ be a normal projective variety of dimension $3$,
$x_0 \in X$ a nonhypersurface terminal singular point
for ${\ind}_{x_0}X = r \geq 2$, 
and $L$ an ample ${\mathbb Q}$-Cartier divisor such that 
$K_X + L$ is Cartier at $x_0$.
Assume that there are positive numbers ${\sigma}_p$ for $p = 1,2,3$
which satisfy the following conditions: 
$(1)$ $\sqrt[p]{ (L)^p \cdot W} \geq {\sigma}_p$ 
for any subvariety $W$ of dimension $p$ which contains $x_0$,
$(2)$ ${\sigma}_1 \geq 1 + 1 / r$,
${\sigma}_2 \geq (1 + 1 / r){\sqrt{r+3}}$, and 
${\sigma}_3 >  (1 + 1 / r){\sqrt[3]{r+2}}$.
Then $| K_X + L |$ is free at $x_0$.\\
Acknowledgment:
The author would like to express his thanks to
Professor Yujiro Kawamata
for his advice and warm encouragement.
He also would like to express his thanks to
Mr. Masayuki Kawakita 
for teaching him RIMS-1273 [KMMT] preprint's existence.  
\section{Preliminaries}
\begin{defn}  \normalfont
Let ${m_X}_{p}$ be the maximal ideal of $p$ of $X$.
The \textit{embedding dimension} of $X$ at $p$
is the dimension of the Zariski tangent space, 
\[\emb {\dim}_{p} X=
\dim  \frac{ {m_X}_{p} }{ {m_X}^2_{p} }. \]
\end{defn}
We basically use the following: 
\begin{thm}[{[A]}]
Let $S$ be a normal surface, $p$ be a point of $S$.
Suppose $S$ has a rational singularity at $p$.
Let $Z$ be the fundamental cycle.
Then,  
\[{\mult}_{p}S= - Z^2, 
 \textrm{ for all integers $k$  } 
\dim \frac{ {m_S}_{p}^k }{ {m_S}_{p}^{k+1} }
= k {\mult}_{p} S + 1, \]  
\[ \textrm{and }   
 \emb {\dim}_p S = {\mult}_{p} S + 1 .\]
\end{thm}
We would like to calculate the multiplicity of terminal singularities on
threefolds. We shall need the following Mori's classification 
theorem of terminal singularities in dimension $3$. 
\begin{thm}[{[M]}]
Let $0 \in X$ be a $3$-fold terminal nonhypersurface singular point  
over $\mathbb{C}$.
Then $0 \in X$ is isomorphic to a singularity 
described by the following list:\\
$(1)cA/r, 
\{xy+f(z,u^r)=0, f \in \mathbb{C}\{ z, u^r \}, (r, a)=1 \}
\subset \mathbb{C}^4/\mathbb{Z}_r(a,r-a,r,1)$, \\
$(2)cAx/4, \{ x^2+y^2+f(z, u^2)=0, f \in \mathbb{C}\{ z,u^2 \} \}
\subset \mathbb{C}^4/\mathbb{Z}_4(1,3,2,1)$, \\
$(3)cAx/2, \{ x^2 + y^2 +f(z, u)=0,  
f \in (z,u)^4\mathbb{C}\{ z, u \} \}
\subset \mathbb{C}^4/\mathbb{Z}_2(1,2,1,1)$, \\
$(4)cD/2, \{u^2+z^3+x y z+f(x, y)=0, f \in (x, y)^4\}, 
\textrm{or } \{u^2+x y z+z^n+f(x,y)=0, f \in (x, y)^4, n \geq 4\}, 
\{{u}^2 + {y}^2 {z} + {z} ^n+f(x, y)=0, f\in(x,y)^4, n \geq 3 \}
\subset \mathbb{C}^4 / \mathbb{Z}_2(1,1,2,1)$, \\
$(5)cD/3, u^2+x^3+y^3+z^3=0, \textrm{or } 
\{u^2+x^3+y z^2+f(x, y, z)=0, f\in (x, y, z)^4\}, 
\{u^2+x^3+y^3+f(x, y, z)=0, f\in (x, y, z)^4 \}
\subset \mathbb{C}^4 / \mathbb{Z}_3(1, 2, 2, 3)$, \\
$(6)cE/2, \{ u^2+x^3 +g(y, z) x+ h(y, z)=0, g,h 
\in \mathbb{C}\{y, z\}, g, h \in (y,z)^4 \}
\subset \mathbb{C}^4 /\mathbb{Z}_2(2,1,1,1)$. \\
The equations have to satisfy $2$ obvious conditions:
$1$. The equations define a terminal hypersurface singularity.
$2$. The equations are ${\mathbb Z}_n$-equivariant. 
$(\textrm{In fact\ }{\mathbb{Z}_n}\textrm{-invariant,\ except\ for\ }cAx/4.)$
\end{thm} 
\section{Main Theorem}
\begin{thm}
Let $(X,p)$ be a $3$-fold terminal singular point 
over $\mathbb{C}$. Then,
for all integers $k$
\[{\dim \  }  \frac{ {m_{X_p}}^k }{ {m_{X_p}}^{k+1}} 
= {\mult}_{p} X \cdot \frac{ k ( k + 1 ) }{2} + k+1, \ 
{\emb}{\dim}_{p} X = {\mult}_{p} X + 2, \]
\[ \textrm{and } {\mult}_{p} X \leq {\ind}_{p} X + 2 \ 
( \textrm{\ if\ }{\ind}_{p} X =1, 
\textrm{\ then\ } {\mult}_{p} X \leq 2). \]
Moreover, we assume that $(X,p) \cong ( x y + f(z, {u}^r) = 0 
\subset {\mathbb{C}}^4/ \mathbb{Z}_r(a, r-a, r, 1), 0)$
or $( \mathbb{C}^3 / \mathbb{Z}_r(a, r-a, 1), 0)$ 
for $(r,a)=1$ and $r > 1$.
Let $r_i := \min \{ r_{i-1} - a_{i-1}, a_{i-1} \} 
(r_0 > r_1 > \cdots > r_n = 1 )$ and 
$a_i = r_{i-1} ( \textrm{\ mod\ } r_i )$ for $r_0 = r$, $a_0 = a$. Then, 
\[ {\mult}_p X = 
{ \frac{}{} }_{\lfloor} { \frac{r_0}{r_1} }_{\rfloor} + 
{ \frac{}{} }_{\lfloor} { \frac{r_1}{r_2} }_{\rfloor} + \cdots +
{ \frac{}{} }_{\lfloor} { \frac{{r}_{n-1}}{{r}_{n} } }_{\rfloor} +2 
\leq r +2 
( = \textrm{if and only if } r_1 = 1 ). \]
In othercases
$(\textrm{$(2)$,$(3)$,$(4)$,$(5)$, or $(6)$ of Theorem 1.3 } )$,
then ${\mult}_p X = r+2$.
\end{thm}
\begin{proof}
Case $0$. 
Let $(X,p)$ be a smooth point. It is clear.\\
Case $1$. 
Let $(X,p)$ be a Gorenstein terminal singular point. 
Since we have $t^2 +f(x,y,z)=0$, 
then $(x,y,z)^k/(x,y,z)^{k+1}$ or 
$t \cdot (x,y,z)^{k-1}/(x,y,z)^{k} \in m^k/m^{k+1}$.
Hence, 
\[ {\dim \ }  \frac{ m_{X_p}^k }{ m_{X_p}^{k+1} } 
= 2 \frac{(k+2)(k+1)}{2} - (k+1) = 2 \frac{k(k+1)}{2} +k+1. \]
\[ \textrm{Hence\ } {\mult}_{p} X = 2 \textrm{\ and\ } 
{\emb}{\dim }_{p} X ={\mult}_{p} X + 2 = 4. \]
Case $2$.
Let $(X,p)$ be a terminal quotient singular point of type \\ 
${\mathbb C}^3/{\mathbb Z}_r(a,-a,1)$ with $(r,a)=1$
Let $S_i={\mathbb C}^2/{\mathbb Z}_r(i,1)$
for $i = a, -a$. For $i=a, -a$, we have 
\[ (xy)^w \cdot \frac{ m_{{S_i}_p}^{k-w} }{ m_{{S_i}_p}^{k-w+1} }
\in \frac{ m_{X_p}^k}{m_{X_p}^{k+1}} \textrm{\ and\ } 
\frac{m_{{S_a}_p}^{k-w}}{ m_{{S_{a}}_p}^{k-w+1} } \cap 
\frac{ m_{{S_{-a}}_p}^{k-w} }{ m_{{S_{-a}}_p}^{k-w+1} } = (z^r)^{k-w} 
\textrm{\ for\ } 0 \leq w \leq k.\] Then by Theorem 1.2, 
\[ {\dim \  } \frac{ {m_X}_p^{k} }{ {m_X}_p^{k+1} }  
= \sum_{w=0}^{k}
 \{ \dim \  \frac{ {m_{S_{a}}}_p^{k-w} }{ {m_{S_{a}}}_p^{k-w+1} }
 +\dim \ \frac{ {m_{S_{-a}}}_p^{k-w} }{ {m_{S_{-a}}}_p^{k-w+1}} -1 \}= \]
\[ \sum_{w=0}^{k}
 \{ (k-w)({\mult}_{p} S_{a} + {\mult}_{p} S_{-a}) + 1 \}
= ({\mult}_{p} S_{a} + {\mult}_{p} S_{-a}) \cdot \frac{k(k+1)}{2} + k + 1. \]  
Hence, 
${\mult}_{p} X = {\mult}_{p} S_{a} + {\mult}_{p} S_{-a}$ 
and ${\emb}{\dim}_{p} X ={\mult}_{p} X +2$. \\
Since we have 
\[ \frac{r_{i}}{r_{i+1}} 
=( { \frac{}{} }_{\lfloor} { \frac{r_{i}}{r_{i+1}} }_{\rfloor} + 1 ) 
- \frac{r_{i+1} - a_{i+1} }{r_{i + 1}}, \textrm{\ then\ } \]
\[ \mult_{p} {\mathbb{C}}^2 / {\mathbb{Z}}_{r_{i}} (r_{i+1},1)
= { \frac{}{}}_{\lfloor} { \frac{r_{i}}{r_{i+1}} }_{\rfloor} - 1 + 
\mult_{p}{\mathbb{C}}^2 / {\mathbb{Z}}_{r_{i+1}} (r_{i+1}-a_{i+1},1). \]
Since we have that
\[ \frac{r_{i} }{r_{i} - r_{i+1}}
=2 - \frac{r_{i} - 2 r_{i+1}}{r_{i} - r_{i+1}}, 
\textrm{\ that\ } \frac{r_{i}-r_{i+1}}{r_{i} - 2 r_{i+1}}
=2 - \frac{r_{i} - 3 r_{i+1}}{r_{i} - 2 r_{i+1}},  
\cdots, \]  
\[ \textrm{\ and\ that\ }
\frac{ r_{i}-( {}_{\lfloor} {{r_{i} }/{r_{i+1} } }_{\rfloor} - 2)r_{i+1}}
{r_{i} - ( {}_{\lfloor} {r_{i}/r_{i+1}}_{\rfloor} -1) r_{i+1}}
=2 - 
\frac{ r_{i} - ( {}_{\lfloor} { r_{i}/r_{i+1}}_{\rfloor} ) r_{i+1} }
{r_{i} -  ( {}_{\lfloor} {r_{i}/r_{i+1}}_{\rfloor} - 1) r_{i+1}},\] 
\[ \textrm{then\ } 
\mult_{p} {\mathbb{C}}^2 / {\mathbb{Z}}_{r_{i}} (r_{i} - r_{i+1},1) 
= 1 + \mult_{p} {\mathbb{C}}^2 / {\mathbb{Z}}_{r_{i+1}}(a_{i+1}, 1). \]
Thus  $\mult_{p} {\mathbb{C}}^2 / {\mathbb{Z}}_{r_{i}} (r_{i+1},1)
+ \mult_{p} {\mathbb{C}}^2 / {\mathbb{Z}}_{r_{i}} (r_{i} - r_{i+1},1) = \\ 
{}_{\lfloor} {r_{i}/r_{i+1}}_{\rfloor} +  
\mult_{p}{\mathbb{C}}^2 / {\mathbb{Z}}_{r_{i+1}} (r_{i+2},1)
+\mult_{p} {\mathbb{C}}^2 / {\mathbb{Z}}_{r_{i+1}}(r_{i+1}-r_{i+2}, 1)$.
Hnece,  
\[ {\mult}_{p} X = { \frac{}{} }_{\lfloor} { \frac{r_0}{r_1} }_{\rfloor} +
 { \frac{}{} }_{\lfloor} { \frac{r_1}{r_2} }_{\rfloor} + \cdots 
+ { \frac{}{} }_{\lfloor} { \frac{r_{n-1}}{r_n} }_{\rfloor} + 2 \leq 
r + 2\ 
(\textrm{\ last\ } = \textrm{\ if\ and\ only\ if\ } r_1 =1 ). \]
Case $3$. 
Let $(X,p)$ be a $3$-fold terminal nonhypersurface singular point. 
We shall use the Mori's classification theorem 
of terminal singularities in dimension $3$ ([M]). \\
Case $3$-$1$. 
$(1)cA/r, \{x y + f(z, {u}^r) = 0, f \in \mathbb{C}\{ z,u^r \} (r,a)=1 \}
\subset \mathbb{C}^4 / \mathbb{Z}_r(a, r - a, r, 1)$. \\
Let $S_i= \mathbb{C}^2/ \mathbb{Z}_r(i/r,1/r)$ for $i=a, -a$.
For $i=a, -a$, we have 
\[z^w \cdot \frac{ m_{{S_i}_p^{k-w}} }{ m_{{S_i}_p^{k-w+1}} }
\in \frac{ m_{X_p}^k }{ m_{X_p}^{k+1} } \textrm{\ and\ } 
\frac{ {m_{S_a}}_p^{k-w} }{ {m_{S_{a}}}_p^{k-w+1} } \cap 
\frac{ {m_{S_{-a}}}_p^{k-w} }{ {m_{S_{-a}}}_p^{k-w+1} } = (u^r)^{k-w}  
\textrm{\ for\ } 0 \leq w \leq k. \]
The rest of the proof is the same as Case $2$.
Hence, 
${\emb}{\dim}_{p} X ={\mult}_{p} X +2$, 
\[ {\mult}_{p} X = 
 {\frac{}{} }_{\lfloor} { \frac{r_0}{r_1} }_{\rfloor} +
 {\frac{}{} }_{\lfloor} { \frac{r_1}{r_2} }_{\rfloor} + \cdots + 
 {\frac{}{} }_{\lfloor} { \frac{r_{n-1}}{r_n} }_{\rfloor} + 2 
\leq {\ind}_p X + 2 
(\textrm{\ last\ } = \textrm{\ if\ and\ only\ if\ } r_1 = 1 ). \]
Case $3$-$2$.
$(2)cAx/4, \{x^2+y^2 +f(z,u^2)=0, f \in \mathbb{C}\{ z,u^2 \} \}
\subset {\mathbb{C}}^4 / {\mathbb{Z}}_4(1, 3, 2, 1)$.\\
We have  $m_{X_p}/m_{X_p}^{2}=(yu,yx,u^4,u^2z,z^2,xu^3,xuz,x^2z)$.
and ${\emb}{\dim}_{p} X =8$.\\
Then, ${(y u)}^t {u}^{4 ( k - t ) - 2 s } {z}^s (0 \leq s \leq 2 ( k - t ) ),
{(y u)}^t x {u}^{4 ( k - t ) - 2 s - 1 } {z}^s 
( 0 \leq s \leq 2 ( k - t ) - 1),
{(y u)}^t ({x}^2 z){(y x)}^s {({z}^2)}^{k - s - t - 1}
(0 \leq s \leq k - t - 1)$, and 
${(y u)}^t {(y x)}^s {({z}^2)}^{k - s - t}(1 \leq s \leq k - t)
 \in m_{X_p}^k / m_{X_p}^{k+1}$.
\[ \dim \ \frac{ {m_X}_p^{k} }{ {m_X}_p^{k+1} }
= \sum_{t=0}^{k} \{(2 k - 2 t + 1) + ( 2 k - 2 t ) 
+ ( k - t ) + ( k - t) \} \]
\[ = \sum_{t=0}^{k} (6 k - 6 t + 1) = 6 \frac{k ( k + 1)}{2} + k +1. \]
Hence, ${\mult}_{p} X = {\ind}_{p}X + 2 =4 + 2 = 6$ 
and ${\emb}{\dim}_{p} X = {\mult}_{p} X + 2 = 8$. \\
Case $3$-$3$. $(3)(4)(6)$ \\
$(3)cAx/2, \{{x}^2 + {y}^2 + f(z, u) = 0, 
f \in (z,u)^4\mathbb{C}\{ z,u \} \} \subset 
{\mathbb{C}}^4 / {\mathbb{Z}}_2(1, 2, 1, 1)$.\\
We have $m_{X_p}/m_{X_p}^{2}=(y,z^2,zu,u^2,xz,xu)$ and 
${\emb}{\dim}_{p} X =6$.\\
Then, ${y}^s (z, u)^{2 (k - s)} (0 \leq s \leq k),
x {y}^t (z, u)^{2 (k - t) - 1}(0 \leq t \leq k - 1)
\in m_{X_p}^k/m_{X_p}^{k+1}$.
\[ {\dim \  } \frac{ {m_X}_p^{k} }{ {m_X}_p^{k+1} }
= \sum_{s=0}^{k}(2k-2s +1) +\sum_{t=0}^{k-1}(2k-2t)  
= 4 \frac{k ( k + 1)}{2}+ k +1. \]
Hence ${\mult}_{p} X ={\ind}_{p} X + 2 = 4$ 
and ${\emb}{\dim}_{p} X = {\mult}_{p}X + 2 = 6$. \\
The proofs of $(4)$ and $(6)$ are the same 
as the proof of $(3)$.\\
Case $3$-$4$ 
$(5)cD/3, {u}^2 + {x}^3 + {y}^3 + {z}^3 = 0$, or
$\{ {u}^2 + {x}^3 + y {z}^2 + f(x, y, z) = 0, f \in (x, y, z)^4 \}$, 
$\{ {u}^2 + {x}^3 + {y}^3 + f(x, y, z), f\in (x,y,z)^4 \} 
\subset {\mathbb{C}}^4 / {\mathbb{Z}}_3(1, 2, 2, 3)$.\\
We have $m_{X_p}/m_{X_p}^{2}=(xy,xz,u,y^3,y^2z,yz^2,z^3)$.
and ${\emb}{\dim}_{p} X =7$.\\
Then, ${u}^t {(x (y, z))}^{ k - t } (0 \leq t \leq k),
{u}^t {(x (y, z))}^{k - t - 1} {(y, z)}^3 (0 \leq t \leq k - 1 ), \\
{u}^t {(x (y, z))}^s {(y, z)}^{ 6 + 3 ( k - 2 - s - t) }
(0 \leq s + t \leq k - 2 ) \in m_{X_p}^k/m_{X_p}^{k+1}$. 
\[ {\dim \  } \frac{ {m_X}_p^{k} }{ {m_X}_p^{k+1} }=
\sum_{t=0}^{k}(k-t+1) +\sum_{t=0}^{k-1}(k-t+3)+
 \sum_{t=0}^{k-2}\sum_{s=0}^{k-t-2} 3 \]
\[= \frac{1}{2} \cdot (k + 1)(k + 2)+ \frac{1}{2} \cdot k (k+7) +
\frac{3}{2} k (k-1) = 5 \frac{k(k+1)}{2} + k + 1. \]
Hence 
${\mult}_{p}X={\ind}_{p} X + 2 = 5$ and 
${\emb}{\dim}_{p} X ={\mult}_{p} X + 2 = 7$. 
\end{proof}
We give the following concrete example: 
\begin{Exam}
Let $(X,p)$ be a quotient singular point of type
${\mathbb C}^3/{\mathbb Z}_{13}(5,8,1)$.
Then, $\mult_{p} X = {}_[ (13/5)_] + {}_[ (5/2)_] + {}_[ (2/1)_] + 2 = 8$.
\end{Exam}
Theorem 2.1 is wrong on the following canonical singularity on threefolds.  
\begin{Exam}
Let $(X,p)$ be a quotient singular point of type
${\mathbb C}^3/{\mathbb Z}_3(1,1,1)$.
We have ${\mult}_{p} X =9$ and ${\emb}{\dim}_{p} X=10$
Then, ${\emb}{\dim}_{p}X = 10 < {\mult}_{p} X + 2 = 11$ 
and ${\mult}_{p} X = 9 > {\ind}_{p} X + 2 = 5$.
\end{Exam}
\section{Application1}
We can improve the boundedness number 
in [KMMT,Theorem 1.2 (2)] by Theorem 2.1.  
\begin{defn}[KMMT Theorem 1.2]\normalfont
Let $X$ be a normal projective variety and
$X$ is called a \textit{terminal} $( \textrm{resp.} \textit{klt} )$ 
\textrm{$\mathbb{Q}$-Fano variety},
if $X$ has only terminal singularities and $-K_X$ is ample.
By replacing 'ample' with 'nef and big', 
\textit{terminal} $( \textrm{resp.} \textit{klt} )$ 
\textit{weak $\mathbb{Q}$-Fano varieties} are similarly defined.
Let $I(X)$ be the smallest  positive integer $I$ 
such that $I K_X$ is Cartier;
$I(X)$ is called \textit{the Gorenstein index} of $X$.
We note that if $X$ is a klt $\mathbb{Q}$-Fano variety then $|-m K_X|$ is
free for some $m >0$.
The induced birational morphism $X \rightarrow \bar{X}$
is said to be the \textit{anti-canonical morphism} of $X$.
\end{defn}
\begin{Lemm}[{[KMMT Lemma 4.1]}]
Let $X$ be an $n$-dimensional projective variety and $x$ a closed
point with multiplicity $r$.
Let $D$ be a nef and big $\mathbb{Q}$-Cartier divisor on $X$
and $l$ a covering family of curves containing $x$ 
such that $D \cdot l \leq d$. Then $D^n \leq rd^n$.
\end{Lemm}
The following is our improvement for [KMMT Theorem 5.1].
\begin{thm}
Let $X$ be a $\mathbb{Q}$-factorial terminal 
$\mathbb{Q}$-Fano $3$-fold with $\rho (X)=1$.
Then $(-K_X)^3 \leq 6^3 \cdot \ ( 2 + 24!)$.
\end{thm}
\begin{proof}
(cf. [KMMT Theorem 5.1])
By [MM 86, Thm.5], there is a covering family of rational curves
$\{l\}$ such that $-K_X \cdot l \leq 6$.
If $\{l \}$ has a fixed point $x$, then Lemma 3.2,
we have $(-K_X)^3 \leq 6^3 \cdot {\mult}_{x}X$.\\
We have ${\mult}_{x}X  \leq  2 + {\ind}_{x} X$.
By [KMMT Theorem 1.2 (1)], 
we have ${\ind}_x X \leq 24!$.
Hence $(-K_X)^3 \leq 6^3 \cdot (2 + 24!)$ in this case.

If $\{l \}$ has a fixed point $x$,
the proof is the same as the one of [KMM92a, Theorem.].

By [KMMT, Construction-Proposition 4.4 and Claim 5.2], 
there is a covering family of rational curves $\{ l' \}$ 
with a fixed point $x$ such that $-K_X \cdot \{ l' \} \leq 3 \times 6$. 
Hence by Lemma 3.2, $(-K_X)^3 \leq 6^3 \cdot 3^3$ in this case.
\end{proof}
The following is our improvement for [KMMT Theorem 1.2].
\begin{thm}
Let $X$ be a terminal weak $\mathbb{Q}$-Fano $3$-fold.
Then the following hold.
$(1)$ $ -K_x \cdot c_2(X) \geq 0$, and hence $I(X)|24!$. 
$(2)$ Assume further that the anti-canonical morphism 
$g:X \rightarrow  \bar{X}$ does not contract any divisors. 
Then $(-K_X)^3 \leq 6^3 \cdot (2 + 24!)$.
$(3)$The terminal $\mathbb{Q}$-Fano $3$-folds are bounded.
\end{thm}
\begin{proof}
The proof is the same as the one of [KMMT Theorem 1.2]
except that 
we can use Theorem 3.3 instead of [KMMT Theorem 5.1].
\end{proof}
\section{Application2}
We can extend [K,Theorem 3.6 ] to nonhypersurface terminal
singularities in the following. 
\begin{thm}
Let $X$ be a normal projective variety of dimension $3$,
$x_0 \in X$ a nonhypersurface terminal singular point
for ${\ind}_{x_0}X = r \geq 2$, 
and $L$ an ample ${\mathbb Q}$-Cartier divisor such that 
$K_X + L$ is Cartier at $x_0$.
Assume that there are positive numbers ${\sigma}_p$ for $p = 1,2,3$
which satisfy the following conditions: \\
$(1)$ $\sqrt[p]{ (L)^p \cdot W} \geq {\sigma}_p$ 
for any subvariety $W$ of dimension $p$ which contains $x_0$,\\
$(2)$ ${\sigma}_1 \geq 1 + 1 / r$,
${\sigma}_2 \geq (1 + 1 / r){\sqrt{ r + 3 }}$, and 
${\sigma}_3 >  (1 + 1 / r){\sqrt[3]{ r + 2 }}$.\\
Then $| K_X + L |$ is free at $x_0$.
\end{thm}
\begin{proof}
We have ${\mult}_{x_0} X \leq r+2$ and 
${\emb}{\dim}_{x_0} X  \leq r+4$.
The rest of the proof is the same as the one of [K, Theorem 3.6.]. 
\end{proof}

\end{document}